\documentclass[11pt]{article}
\usepackage{graphicx}
\usepackage{amsmath}
\usepackage{amssymb}
\usepackage{amsthm}
\usepackage[titletoc]{appendix}
\usepackage{url}
\usepackage{enumerate}
\newcommand\RR{\mathbb{R}}

\newcommand\al\alpha
\newcommand\be\beta
\newcommand\de\delta
\newcommand\ep\varepsilon
\newcommand\tha\theta
\newcommand\ka\kappa
\newcommand\la\lambda
\newcommand\om\omega

\newcommand\iy\infty
\newcommand\pa\partial

\newcommand{\hyp}[5]{\,\mbox{}_{#1}F_{#2}\!\left(\genfrac{}{}{0pt}{}{#3}{#4};#5\right)}

\numberwithin{equation}{section}
\newtheorem{theorem}{Theorem}[section]

\newtheorem{Remark}[theorem]{Remark}
\newenvironment{remark}{\begin{Remark}\rm}{\end{Remark}}

\begin{document}

\title{An extension of the orthogonal derivative with adjustable precision}
\author{Enno Diekema \footnote{email adress: e.diekema@gmail.com}}
\maketitle

\begin{abstract}
\noindent
The orthogonal derivative is defined as a limit of an integral whose kernel contains an orthogonal polynomial with its measure. When in practice no limit is taken, it means that the accuracy of the derivative depends on the second derivative of the given function. Liptaj shows that it is possible to define a kernel in such a way that the accuracy depends on a higher derivative at your own choice. The accuracy is therefore much greater than with the orthogonal derivative. However Diekema and Koornwinder find a similar extension starting directly from of the orthogonal derivative. The new kernel is not orthogonal for order greater then one. The transfer function for this new derivative is given.
\end{abstract}

\section{Introduction}
\setlength{\parindent}{0cm}
In many applications one needs to approximate the first or higher order derivative of a function which is perturbed by noise. A good candidate for an approximation of the first derivative $f'(x)$ of a continuous function $f(x)$ is given by
\[
f'(x) \approx -\dfrac{3}{2\delta}\int_{-1}^1f(x+\delta\xi)\xi d\xi
\]
for small $\delta$. So part of the noise will be averaged out, because the integral has an averaging effect. The history of this formula is described in \cite{2}. One can recognise that the function $\xi$ is the first Legendre polynomial.
In the paper \cite{2} this idea is generalized to the so-called orthogonal derivative. For the first order the derivative is called the Lanczos' derivative.  

In this paper a modification of the orthogonal derivative is described. This modification has the great advantage that the accuracy of the approximate derivative can be set by means of an extra parameter. There are two methods to go about this modification.  The first method is based on an idea in some papers of Liptaj \cite{3,4}. This method is based on the basic idea of integration by differentiation, which originated from the method of the least squares. The second method uses a formula from \cite[Section 4]{2}. This formula originates from the use of orthogonal polynomials, which is also essentially based on the use of the least squares method. Both methods derive a formula for the kernel of an integral. The first method is more complicated than the second method.

The content of this paper is as follows. The orthogonal derivative is discussed in section 2. This is a part of the first chapter from the author's thesis \cite{1}. Section 3 derives the general formula for the kernel of the approximate derivative of arbitrary order following the basis idea of Liptaj. The result is a formula with a summation of a hypergeometric function. The final polynomial is not orthogonal for the orders greater than one. So the derivative cannot be an orthogonal derivative. In section 4 some special cases for the basic kernel are described. In section 5 the second method using a formula from \cite[Section 4]{2} is treated. Section 6 defines the general kernel formula for the first derivative. It turns out that this kernel is an orthogonal polynomial. The known properties of this polynomial are indicated. The formula for the first derivative can be mentioned as an orthogonal derivative or Lanczos' derivative. In section 7 the transfer function of the first derivative is determined. The transfer function is illustrated by a log-log plot. A summary and conclusion is presented in section 8. In the appendices an integral of a Gegenbauer polynomial is derived and an overview is given for the basic kernel function for different values of the order and the degree.

\section{The orthogonal derivative}
\setlength{\parindent}{0cm}
A well-known problem in signal analysis is to detect the signal that was disturbed by noise. In his practice of signal analysis the author met the problem of differentiating a signal that was disturbed by noise. Because differentiating a signal will amplify the higher frequencies in the signal more than the lower frequencies, the result of differentiating a signal with noise becomes unstable. Therefore, an integrating factor must be added. Every averaging has an integrating effect. The author choose in his work for the method of the least square. (see Figure 1).
\begin{figure}[ht]
	\centering     
	\includegraphics[width=6cm]{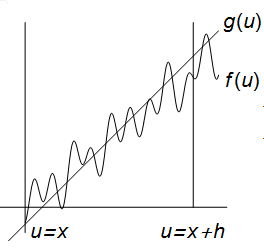} 
	\caption{The definition of the least square derivative.}
	\label{Figure 1.3}
\end{figure}

We approximate the function $f(x)$ by a linear function $g(x)$ with the property that
\begin{equation}
	\int_x^{x+h}\big[f(u)-g(u)\big]^2du
	\label{I.0}
\end{equation}
is minimal. This is the so-called method of the least squares. From Figure \ref{Figure 1.3} it is seen that the noise will be averaged out. We take $g'(x)$ as an approximation of $f'(x)$. With $g(u)=a+bu$ we have to minimize
\begin{equation}
\int_0^1\big[f(x+uh)-a-b(x+uh)\big]^2du
\label{I.1}
\end{equation}
as a function of $a$ and $b$. This can be handled by putting the derivatives of \eqref{I.1} with respect to $a$ and $b$ equal to zero and solving the resulting system of linear equations in $a$ and $b$. The solution for $b$ is
\[
b=\dfrac{6}{h}\int_0^1f(x+hu)(2u-1)du.
\]
Then take $b$ as the approximation of $f'(x)$. If we want to minimize \eqref{I.0} for $g$ a quadratic function then we have to solve a system of three linear equations. More generally, in order to minimize
\eqref{I.0} for $g$ a polynomial of degree $n$ we have to solve a system of $n+1$ linear equations. When doing these computations for some successive low values of $n$ one will see a structure arising in the expression for the approximation of the $n$-th derivative. It involves an integral for which the integrand contains a so-called shifted Legendre polynomial $P_n^*(x)$. Such polynomials are orthogonal polynomials on the interval $[0,1]$ with respect to a constant weight function. The resulting formula for the $n$-th derivative is
\[
\dfrac{d^n}{dx^n} f(x) =\dfrac{\Gamma (2n+2)}{\Gamma (n+1)}\lim\limits_{h\downarrow 0}\dfrac{1}{h^{n}}
\int\nolimits_{x}^{x+h}f(u) P_{n}^{*}\left(\dfrac{u-x}{h} \right) du.
\]
Of course, this formula, obtained for low $n$ by tedious computations, cannot be an accident. Indeed, a simple argument which involves the Taylor series expansion of $f(u+h)$ up to the term with $h^n$ and some very basic properties of orthogonal polynomials, gives a proof of the formula for general $n$. Moreover, the proof of this formula does not depend on the special choice of interval and weight function (or measure). We may start with minimizing
\begin{equation}
	\int_{\RR}\big[f(x+hu)-g(u)\big]^2 d\mu(u)  
	\label{I.2}
\end{equation}
for $g$ a polynomial of degree $\leq n$. $\mu$ is a positive measure on $\RR$, for instance $d\mu(u)=w(u) du$ on some interval and $=0$ outside that interval, where $w$ is a positive weight function. Then it can be shown \cite[Ch. 2.3]{1}, \cite{2} that
\begin{equation}
\dfrac{d^n}{dx^n} f(x) =\lim\limits_{h\downarrow 0}\dfrac{k_{n}n!}{h_{n}}
\dfrac{1}{h^{n}}\int_{\RR}f(x+hu) p_{n}(u) d\mu (u)
\label{I.3}
\end{equation}
where the polynomials $p_n$ are orthogonal polynomials with respect to the measure $\mu$ i.e, $p_n$ is a polynomial of degree $n$ with coefficient of highest degree term equal to $k_n$, and
\[
\int_{\RR} p_m(x) p_n(x) d\mu(x) = h_n \delta_{m,n}.
\]
We call the right-hand side of \eqref{I.3}, if the limit exists, the orthogonal derivative of order $n$ of $f$ at $x$. This limit certainly exists if $f$ is $n$ times differentiable at $x$, and then the orthogonal derivative equals the ordinary derivative, but the limit exists much more generally, for instance, it is equal to the $n$-th order Peano derivative of $f$ at $x$, if this exists.
The following items are some examples of orthogonal derivatives with integration interval $[-1,1]$.
\begin{itemize}
\item Legendre derivative
\begin{equation}
\dfrac{d^n}{dx^n} f(x) =\lim\limits_{h\downarrow 0}\dfrac{k_{n}n!}{h_{n}}
\dfrac{1}{h^{n}}\int_{-1}^1 P_{n}(u)f(x+hu)du.
\label{I.5}
\end{equation}
$P_n(x)$ are the Legendre polynomials with weight function $1$ on the interval $[-1,1]$.
\item Gegenbauer derivative
\[
\dfrac{d^n}{dx^n} f(x) =\lim\limits_{h\downarrow 0}\dfrac{k_{n}n!}{h_{n}}
\dfrac{1}{h^{n}}\int_{-1}^1(1-u^2)^{\alpha-1/2}C^{(\alpha)}_n(u)f(x+hu) du.
\]
$C^{(\alpha)}_n(x)$ are the Gegenbauer polynomials with weight function $(1-x^2)^{\alpha-1/2}$ on the interval $[-1,1]$.
\item Jacobi derivative
\[
\dfrac{d^n}{dx^n} f(x) =\lim\limits_{h\downarrow 0}\dfrac{k_{n}n!}{h_{n}}
\dfrac{1}{h^{n}}\int_{-1}^1(1-u)^\alpha(1+u)^\beta P^{(\alpha,\beta)}_n(u)f(x+hu) du.
\]
$P^{(\alpha,\beta)}_n(x)$ are the Jacobi polynomials with weight function $(1-x)^\alpha(1+x)^\beta$ on the interval $[-1,1]$.
\end{itemize}

As an approximation of the orthogonal derivative we define the $\text{\em approximate orthogonal derivative}$ where we omit the limit for $h \downarrow 0$. This has the property that the noise of a signal will be suppressed. We can illustrate this property by working in the frequency domain. We look at the approximate orthogonal derivative as a filter. This filter can be characterized by its transfer function in the frequency domain.

A filter produces from an input signal an output signal. It has the property that the quotient of the Fourier transform of the output signal and the Fourier transform of the input signal is independent of the input. This quotient, denoted by $H(\omega)$ with $\omega$ the frequency, is called the {\em transfer function} of the filter. It is usually a complex-valued function of the frequency. Both its modulus and its argument are of interest. In this paper we use the modulus of the transfer function.

\begin{figure}[ht]
	\centering     
	\includegraphics[width=9cm]{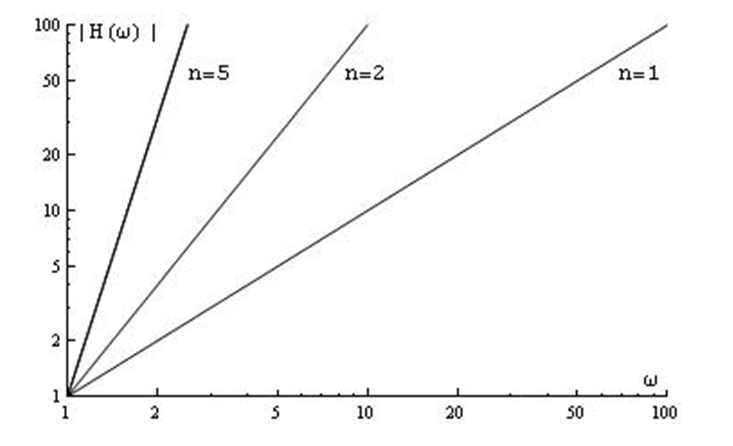} 
	\caption{Moduli of the transfer functions of a differentiator of orders $1$,$2$ and $5$.}
	\label{Figure 1.5}
\end{figure}

For a filter giving the $n$-th order ordinary derivative of the input we have $H(\omega)=(i\omega)^n$, so for the modulus of the transfer function we can draw the following log-log plots of Figure~ \ref{Figure 1.5}.

We see that for high frequencies the modulus of the transfer function goes to infinity. So the high frequencies will be stronger amplified than the lower frequencies. This high frequency noise will disturb the derivative of the original signal very much.

\begin{figure}[ht]
	\centering     
	\includegraphics[width=9cm]{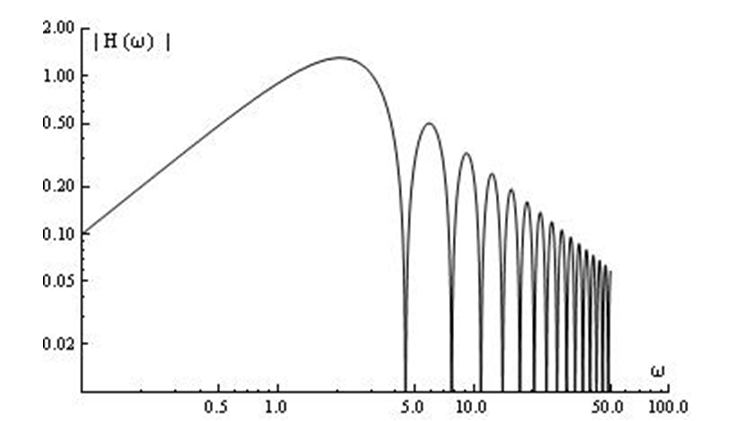} 
	\caption{Modulus of the transfer function of a first order Legendre differentiator.}
	\label{Figure 1.6}
\end{figure}

However, if we use for $n=1$ a filter corresponding to the approximation of the first derivative implied by the right-hand side of \eqref{I.3} with $d \mu(u)=du$ on $[-1,1]$ (involving the first degree Legendre polynomial $P_1(x)=x$) then we have a filter sending input function $f$ to output function $g$ by
\[
g(x) =\dfrac{3}{2h }\int\nolimits_{-1}^{1}f(x+hu) u\ du.
\]
The modulus of the corresponding transfer function is drawn in the log-log plot of Figure~\ref{Figure 1.6}.

\ 

We see that $|H(\omega)| \approx \omega$ for lower frequencies but that it has successive oscillations decreasing to zero for higher frequencies. Therefore, higher frequency parts of the signal are suppressed by the filter.

\section{First method for an alternative for the orthogonal derivative}
Liptaj in a number of papers \cite{3},\cite{4} developed a method for the derivative by integration of a function on an bounded interval with greater accuracy then the approximate orthogonal derivative. 

Looking for example to formula \eqref{I.5} we see that when expanding the function $f(x+h\ t)$ into a series of powers of $h$, the first term we may neglect is of order $n+1$ for a derivative of order $n$. Liptaj found that it is possible to choose a function $k(t)$ in such a manner that the first term to neglect is of higher order than $n+1$. He defined the following functions:
\[
k(t)=\dfrac{d^n}{dt^n}\omega(t).
\]
\[
\omega(t)=K(1-t^2)^n\sum^m_{k=0}a_{2k}\, t^{2k}.
\]
He choose $n$ for the order of the derivative. His formula for the $n$-th order derivative becomes
\begin{equation}
\dfrac{d^n}{dx^n} f(x) =\lim_{h \rightarrow 0}\left(-\dfrac{1}{h}\right)^n\int_{-1}^1 k(t)f(x+h\,t)dt.
\label{2.1}
\end{equation}

\

We illustrate his method by the following example.
We want to find the formula for the third order derivative of the function $f(x)=x^9$. The exact answer is known: $\dfrac{d^3}{dx^3}x^9=7*8*9*x^6$. For the given function we use: $f(x+h\, t)=(x+h\,t)^9$. Theoretically, $h$ should go to zero. In practice, an approximation is used. For the approximate derivative, the following formula is used:
\[
k(t)=K\dfrac{d^3}{dt^3}\Big((1-t^2)^3(1+a_2t^2+a_4t^4+a_6t^6)\Big).
\]
It turns out that with this kernel, the first term to neglect is of order $h^8$. We can use more or less terms, but there are only terms with even power. Using \eqref{2.1} gives
\begin{equation}
\dfrac{d^3}{dx^3}(x+h\,t)^9=K\lim_{h \rightarrow 0}
\left(-\dfrac{1}{h}\right)^n\int_{-1}^1
\dfrac{d^3}{dt^3}\Big((1-t^2)^3(1+a_2t^2+a_4t^4+a_6t^6)\Big)(x+h\,t)^9dt.
\label{2.0}
\end{equation}
Omitting the limit and working out the integral we get
\begin{align*}
\dfrac{d^3}{dx^3}(x+h\,t)^9
&=K\dfrac{256}{230945}\big(4845+2261a_2+1197a_4+693a_6\big)h^6+ \\
&+K\dfrac{768}{2431}\big(663+255a_2+119a_4+63a_6\big)h^4x^2+ \\
&+K\dfrac{768}{143}\big(143+39a_2+15a_4+7a_6\big)h^2x^4+ \\
&+K\dfrac{256}{715}\big(1287+143a_2+39a_4+15a_6\big)x^6.
\end{align*}
If the first three terms are equal to zero and the factor for the last term is equal to $7 * 8 * 9$, then the parameters $a_2,\, a_4,\, a_6$  and $K$ can be determined. In that case, the inaccuracy of the derivative is proportional to $h^8$. The equations are:
\begin{align*}
&4845+2261a_2+1197a_4+693a_6=0. \\
&663+255a_2+119a_4+63a_6=0. \\
&143+39a_2+15a_4+7a_6=0. \\
&K\big(1287+143a_2+39a_4+15a_6\big)=7*8*9.
\end{align*}
It is evident  that if we choose more unknown $a_{2k}$ we get more equations to solve.
The solution of this system of equations is: 
\[
a_2=-15	\qquad a_4=51 \qquad a_6=-\dfrac{323}{7} \qquad K=-\dfrac{105105}{32768}.
\]
Substitution in $\omega(t)$ gives
\[
\omega(t)=\dfrac{105105}{32768}(1-t^2)^3\left(1-15t^2+51t^4-\dfrac{323}{7}t^6\right).
\]
Taking the third derivative results in
\[
k(t)=\dfrac{45045}{4096}\big(17765t^9-39780t^7+30030t^5-8580t^3+693t \big).
\]
We can use this function in formula \eqref{2.1} and get the approximate third derivative.

\

For each order of the derivative (which are up to order 4) and the number of terms $k$ (up to $k=5$) Liptay in his paper \cite{4} gives a method for deriving the equations.
For each case these equations should be solved. He does not give a general solution for the function $k(t)$.

\

In the rest of this section a general formula for the function $k(t)$ is derived. It is assumed that the given function is analytical in the given domain.  We call $k(t)$ the {\em basic kernel function}. Solving $m$ equations with $m$ unknowns is then avoided. We start with
\begin{align}
\dfrac{d^n}{dx^n}f(x)
&=\lim_{h \rightarrow 0} \left(-\dfrac{1}{h}\right)^n \int_{-1}^1 f(x+h\,t)
\dfrac{d^n}{dt^n}\left((1-t^2)^n K \sum_{k=0}^m a_{2k} t^{2k}\right)dt \nonumber \\
&=\lim_{h \rightarrow 0} \left(-\dfrac{1}{h}\right)^n \int_{-1}^1 f(x+h\,t) K \sum_{k=0}^m a_{2k}
\dfrac{d^n}{dt^n}\left((1-t^2)^n t^{2k}\right) dt.
\label{2.2}
 \end{align}
So we have $m+1$ terms in the summation. See \eqref{2.0}. Because $a_0=1$ there are $m$ unknown parameters to detect. The derivative can be taken by means of Leibnitz's rule.
\[
\dfrac{d^n}{dt^n}\big(f(t)g(t)\big)=\sum_{j=0}^n\binom{n}{j}D^{n-j}[f(t)]D^j[g(t)].
\]
Application of \cite[8.934(1)]{7} gives
\[
D^j\big[(1-t^2)^n\big]=(-2)^j\ j!
\dfrac{(-n)_j}{(j-2n)_j}(1-t^2)^{n-j}C^{(n-j+1/2)}_j(t).
\]
$C^{(n-j+1/2)}_j(t)$ are the Gegenbauer polynomials. $(a)_j$ is the well-known Pochhammer symbol.
\[
D^{n-j}\big[t^{2k}\big]=\dfrac{\Gamma(2k+1)}{\Gamma(2k+1-n+j)}t^{2k-n+j}.
\]
Substitution in \eqref{2.2} gives
\begin{align}
k(t)&= K \sum_{k=0}^m a_{2k}\dfrac{d^n}{dt^n}\left((1-t^2)^n t^{2k}\right)\nonumber \\
&=K \sum_{k=0}^m a_{2k}\dfrac{\Gamma(2k+1)}{\Gamma(2k+1-n)}\sum_{j=0}^n\binom{n}{j}
\dfrac{(-n)_j(-2)^j\ j!}{(2k+1-n)_j(j-2n)_j}t^{2k-n+j}(1-t^2)^{n-j}C^{(n-j+1/2)}_j(t).
\label{2.3}
\end{align}
To find formulae for $a_{2k}$ and $K$ we can write the function $f(x+h\, t)$  as a Taylor series. Then there arise all terms with $(x+h\, t)^\lambda$. We only have to take the terms with $\lambda=2m+n$. Because the function $k(t)$ contains only terms with $\lambda \leq 2m+n$ we have only use the factor $(x+h\, t)^{2m+n}$. The factor  $(x+h\, t)^{2m+n}$ can be written as
\[
(x+h\, t)^{m+n}=\sum^{2m+n}_{i=0}\binom{2m+n}{i}x^{2m+n-i}(h\, t)^i.
\]
To compute the value of $a_{2k}$ it is sufficient to take only the term with  $x^{2k+n-i}$.
Then we define with \eqref{2.1} and \eqref{2.3} the function $A(x)$ as 
\begin{multline}
A(x)= (-1)^n K\sum^{2m+n}_{i=0}\binom{2m+n}{i}
\sum_{k=0}^m a_{2k}\dfrac{\Gamma(2k+1)}{\Gamma(2k+1-n)} \\
\sum_{j=0}^n\binom{n}{j}\dfrac{(-n)_j(-2)^j\ j!}
{(2k+1-n)_j(j-2n)_j}x^{2k+n-i}h^{i-n} 
\int_{-1}^1 t^{i+2k-n+j}(1-t^2)^{n-j}C^{(n-j+1/2)}_j(t)dt.
\label{2.5}
\end{multline}
The integral is derived in Appendix A. Using the result and after some manipulations with the Pochhammer symbols and the Gamma functions we obtain
\begin{multline*}
A(x)=(-1)^n K\dfrac{\Gamma(n+1)}{2}\sum^{2m+n}_{i=0}\binom{2m+n}{i}\big[1+(-1)^{i-n}\big]
\sum^m_{k=0}a_{2k}x^{2k+n-i}h^{i-n} \\
\dfrac{\Gamma\left(\dfrac{i+2k-n+1}{2}\right)}{\Gamma\left(\dfrac{i+2k+n+3}{2}\right)}
\dfrac{\Gamma(2k+1)}{\Gamma(2k+1-n)}
\sum^n_{j=0}\dfrac{(-n)_j(i+2k-n+1)_j}{(2k+1-n)_j}\dfrac{1}{j!}.
\end{multline*}
The last summation is well-known.
\[
\hyp21{-n,b}{c}{1}=\dfrac{(c-b)_n}{(c)_n}.
\]
Application gives
\begin{multline*}
A(x)=(-1)^n K\dfrac{\Gamma(n+1)}{2}\sum^{2m+n}_{i=0}\binom{2m+n}{i}\big[1+(-1)^{i-n}\big]
\sum^m_{k=0}a_{2k}x^{2k+n-i}h^{i-n} \\
\dfrac{\Gamma\left(\dfrac{i+2k-n+1}{2}\right)}{\Gamma\left(\dfrac{i+2k+n+3}{2}\right)}
(-i)_n.
\end{multline*}
For the last Pochhammer symbol we use
\[
(-i)_n=(-1)^n\dfrac{\Gamma(i+1)}{\Gamma(i-n+1)}.
\]
This factor gives only a contribution if $i \geq n$. Application gives
\begin{multline*}
A(x)=K\dfrac{\Gamma(n+1)}{2}\sum^{2m+n}_{i=n}\binom{2m+n}{i}\big[1+(-1)^{i-n}\big]
\sum^m_{k=0}a_{2k}x^{2k+n-i}h^{i-n} \\
\dfrac{\Gamma\left(\dfrac{i+2k-n+1}{2}\right)}{\Gamma\left(\dfrac{i+2k+n+3}{2}\right)}
\dfrac{\Gamma(i+1)}{\Gamma(i-n+1)}.
\end{multline*}
Replacing $i$ with $i+n$ and interchanging the summations (which is allowed because all summations are bounded) gives
\begin{multline*}
A(x)=K\dfrac{\Gamma(n+1)}{2}\dfrac{\Gamma(2m+n+1)}{\Gamma(2m+1)}
\sum^m_{k=0}a_{2k}x^{2k}\dfrac{\Gamma\left(\dfrac{1}{2}+k\right)}
{\Gamma\left(\dfrac{3}{2}+n+k\right)} \\
\sum^{m}_{i=0}\big[1+(-1)^{i}\big]\dfrac{(-2m)_i\left(\dfrac{1}{2}+k\right)_{i/2}}{\left(\dfrac{3}{2}+n+k\right)_{i/2}}
\dfrac{1}{i!}\left(\dfrac{h}{x}\right)^i.
\end{multline*}
The second summation gives a hypergeometric function
\[
\sum^{m}_{i=0}\big[1+(-1)^{i}\big]\dfrac{(-2m)_i\left(\dfrac{1}{2}+k\right)_{i/2}}{\left(\dfrac{3}{2}+n+k\right)_{i/2}}
\dfrac{1}{i!}\left(\dfrac{h}{x}\right)^i=
2\, \hyp32{-m,\dfrac{1}{2}-m,\dfrac{1}{2}+k}{\dfrac{1}{2},\dfrac{3}{2}+n+k}
{\dfrac{h^2}{x^2}}.
\]
For $A(x)$ there rests
\[
A(x)=K\Gamma(n+1)\dfrac{\Gamma(2m+n+1)}{\Gamma(2m+1)}
\sum^m_{k=0}a_{2k}x^{2k}\dfrac{\Gamma\left(\dfrac{1}{2}+k\right)}
{\Gamma\left(\dfrac{3}{2}+n+k\right)}
\hyp32{-m,\dfrac{1}{2}-m,\dfrac{1}{2}+k}{\dfrac{1}{2},\dfrac{3}{2}+n+k}
{\dfrac{h^2}{x^2}}.
\]
This is a polynomial in $(x\, h)^2$. The factors before the powers of $(x\, h)^2$ should be zero. So
\begin{align*}
S&=\sum^m_{k=0}a_{2k}\dfrac{\Gamma\left(\dfrac{1}{2}+k\right)}
{\Gamma\left(\dfrac{3}{2}+n+k\right)}
\hyp32{-m,\dfrac{1}{2}-m,\dfrac{1}{2}+k}{\dfrac{1}{2},\dfrac{3}{2}+n+k}
{\dfrac{h^2}{x^2}} \\
&=\sum^m_{k=0}a_{2k}\dfrac{\Gamma\left(\dfrac{1}{2}+k\right)}
{\Gamma\left(\dfrac{3}{2}+n+k\right)}
\sum_{i=0}^m
\dfrac{(-m)_i\left(\dfrac{1}{2}-m\right)_i\left(\dfrac{1}{2}+k\right)_i}
{\left(\dfrac{1}{2}\right)_i\left(\dfrac{3}{2}+n+k\right)_i}\dfrac{1}{i!} \\
&=\dfrac{\Gamma\left(\dfrac{1}{2}\right)}{\Gamma\left(\dfrac{3}{2}+n\right)}
\sum_{i=0}^m\dfrac{(-m)_i\left(\dfrac{1}{2}-m\right)_i}{\left(\dfrac{3}{2}+n\right)_i}
\dfrac{1}{i!}\sum_{k=0}^ma_{2k}\dfrac{\left(\dfrac{1}{2}+i\right)_k}{\left(\dfrac{3}{2}+n+i\right)_k}=0
\end{align*}
with $a_0=1$. This gives $m$ equations to compute the $a_{2k}$. But it is easy to check that the solution of these equations is
\begin{equation}
a_{2k}=\dfrac{(-m)_k\left(m+n+\dfrac{3}{2}\right)_k}{\left(\dfrac{3}{2}\right)_k k!}.
\label{2.6}
\end{equation}
Using this formula for $\omega(t)$ together with \cite[8.932.3]{7}  gives the polynomial
\begin{align}
\omega(t)=K(1-t^2)^n\sum_{k=0}^m a_{2k} t^{2k}&
=K(1-t^2)^n\hyp21{-m,m+n+\dfrac{3}{2}}{\dfrac{3}{2}}{t^2} \nonumber \\
&=-K(1-t^2)^n\dfrac{\Gamma\left(-m-n-\dfrac{1}{2}\right)\Gamma(m+1)}
{2\Gamma\left(\dfrac{1}{2}-n\right)}\dfrac{1}{t}C^{n+1/2}_{2m+1}(t).
\label{3.7b}
\end{align}
To compute the factor $K$ we set
\[
K\Gamma(n+1)\sum_{k=0}^m a_{2k}\dfrac{\Gamma\left(\dfrac{1}{2}+k\right)}
{\Gamma\left(\dfrac{3}{2}+n+k\right)}=1.
\]
Substitution of $a_{2k}$ gives
\begin{align*}
K\Gamma(n+1)\sum_{k=0}^m a_{2k}\dfrac{\Gamma\left(\dfrac{1}{2}+k\right)}
{\Gamma\left(\dfrac{3}{2}+n+k\right)}
&=K\Gamma(n+1)\sum_{k=0}^m \dfrac{(-m)_k\left(m+n+\dfrac{3}{2}\right)_k}{\left(\dfrac{3}{2}\right)_k k!}
\dfrac{\Gamma\left(\dfrac{1}{2}+k\right)}{\Gamma\left(\dfrac{3}{2}+n+k\right)}= \\
&=K\dfrac{\Gamma\left(\dfrac{1}{2}\right)\Gamma(n+1)}{\Gamma\left(\dfrac{3}{2}+n\right)}
\hyp32{-m,m+n+\dfrac{3}{2},\dfrac{1}{2}}{n+\dfrac{3}{2},\dfrac{3}{2}}{1}.
\end{align*}
The hypergeometric function is 1-balanced. So the summation is well-known. After lots of manipulation with the Pochhammer symbols and the Gamma functions we get
\begin{equation}
K=\dfrac{2}{\pi}\dfrac{\Gamma\left(m+n+\dfrac{3}{2}\right)\Gamma\left(m+\dfrac{3}{2}\right)}{\Gamma(m+1)\Gamma(m+n+1)}.
\label{2.7}
\end{equation}
Substitution in \eqref{3.7b} gives
\[
\omega(t)=\dfrac{(-1)^m}{\pi}\dfrac{\Gamma\left(m+\dfrac{3}{2}\right)\Gamma\left(n+\dfrac{1}{2}\right)}{\Gamma(m+n+1)}\dfrac{(1-t^2)^n}{t}C_{2m+1}^{(n+1/2)}(t).
\]
For the kernel $k(t)$ we combine \eqref{2.3},\eqref{2.6} and \eqref{2.7}. This results in
\begin{multline}
k(t)=\dfrac{2}{\pi}\dfrac{\Gamma\left(m+n+\dfrac{3}{2}\right)\Gamma\left(m+\dfrac{3}{2}\right)}{\Gamma(m+1)\Gamma(m+n+1)} \sum_{k=0}^m \dfrac{(-m)_k\left(m+n+\dfrac{3}{2}\right)_k}{\left(\dfrac{3}{2}\right)_k }\dfrac{1}{k!}t^{2k} \\
\sum_{j=0}^n(-1)^{n-j}\dfrac{(-2k)_{n-j}(-n)_j(-2n)_j}
{2^j\left(\dfrac{1}{2}-n\right)_j}\left(\dfrac{1-t^2}{t}\right)^{n-j}C^{(n-j+1/2)}_j(t).
\label{2.9}
\end{multline}
Mathematica gives a hypergeometric function for the last summation. We get at last for the basic kernel
\begin{multline}
k(t)=\dfrac{2}{\pi}\dfrac{\Gamma\left(m+n+\dfrac{3}{2}\right)\Gamma\left(m+\dfrac{3}{2}\right)}{\Gamma(m+1)\Gamma(m+n+1)}
\sum_{j=0}^m \dfrac{(-m)_j\left(m+n+\dfrac{3}{2}\right)_j}{\left(\dfrac{3}{2}\right)_j }\dfrac{1}{j!}t^{2j-n} \\
\dfrac{\Gamma(2j+1)}{\Gamma(2j+1-n)}\hyp32{-n,j+1,j+\dfrac{1}{2}}{j+1-\dfrac{n}{2},j+\dfrac{1}{2}-\dfrac{n}{2}}{t^2}.
\label{2.10}
\end{multline}
Attempts to simplify this formula have not yielded any results so far. A much simpler formula is found in Section 5.

 \section{A closer look at the basic kernel function}

In this section some of the features of the basic kernel function are discussed. The basic kernel depends on two parameters. The first parameter is $n$. This is the order of the derivative. The case $n=1$ is so special that is will be discussed in the next section. The second parameter is $m$. This parameter represents the number of terms used to approximate the derivative. This gives a direct measure of the accuracy of the derivative. 

\

The first attribute that stands out,  is that the hypergeometric function is zero-balanced. However, this attribute does not lead to a known function here. Next it is seen that 

\

if $n=even$ the function $k(t)=even \qquad k(t)=k(-t)$.

if $n=odd\ \ $ the function $k(t)=odd \qquad\  k(t)=-k(-t)$.

\

Using \eqref{2.10} for computing the function $k(t)$ there arises a problem. The fraction before the hypergeometric function is not convergent if \  $2k+1-n \leq 0$. But the corresponding lower term in the hypergeometric function becomes also negative. The summation should be split into two summations. The first summation concerns the positive values of the relevant factor and the second sum concerns the negative values of the relevant factor. In that case the limit should be taken. First we rewrite the basic kernel function and get 
\begin{multline}
k(t)=\dfrac{2^{n+1}}{\pi}\dfrac{\Gamma\left(m+n+\dfrac{3}{2}\right)\Gamma\left(m+\dfrac{3}{2}\right)}{\Gamma(m+1)\Gamma(m+n+1)}
\sum_{j=0}^m \dfrac{(-m)_j\left(m+n+\dfrac{3}{2}\right)_j}{\left(\dfrac{3}{2}\right)_j }\dfrac{1}{j!}t^{2j-n} \\
\dfrac{\Gamma(j+1)\Gamma\left(j+\dfrac{1}{2}\right)}
{\Gamma\left(j+\dfrac{1-n}{2}\right)\Gamma\left(j+\dfrac{2-n}{2}\right)}
\hyp32{-n,j+1,j+\dfrac{1}{2}}{j+\dfrac{1 -n}{2},j+\dfrac{2-n}{2}}{t^2}.
\label{44.1}
\end{multline}
We have to distinguish between the case where n is even or n is odd. For the case $n$ even the term $j+\dfrac{2-n}{2}$ is an integer. Then we have to split the summation for two cases: $j+\dfrac{2-n}{2}>0$ and $j+\dfrac{2-n}{2} \leq 0$ dependent of the value of $m$. When $n$ odd the term $j+\dfrac{1-n}{2}$ is an integer. Then we have to split the summation for two cases: $j+\dfrac{1-n}{2}>0$ and $j+\dfrac{1-n}{2} \leq 0$ and also dependent of the value of $m$. In both cases we can use the following formula 
\begin{multline}
\dfrac{1}{\Gamma(-M)} \ _{p+1}F_p
\left(\begin{array}{l}
	a_0,\dots,a_p \\
	-M,b_2,\dots,b_p
\end{array};x\right)= \\
=\dfrac{x^{M+1}(a_0)_{M+1}\dots(a_p)_{M+1}}{\Gamma(M+2)(b_2)_{M+1}\dots(b_p)_{M+1}}
\ _{p+1}F_p
\left(\begin{array}{l}
	a_0+M+1,\dots,a_p+M+1 \\
	M+2,b_2+M+1,\dots,b_p+M+1
\end{array};x\right)
\label{44.2}
\end{multline}
with $M$ a positive integer. For the proof of this formula see \cite[Lemma 2]{10}.
Both cases are not discussed further here and are left to the reader.

\

A very special case occurs when $m=0$. We distinguish again the cases with $n$ even and $n$ odd. Starting with $n$ even, equation \eqref{44.1} becomes
\[
k(t)=\dfrac{\Gamma\left(n+\dfrac{3}{2}\right)}{\Gamma(n+1)\Gamma\left(\dfrac{1-n}{2}\right)}\dfrac{1}{t^n}\lim_{m \rightarrow 0}\dfrac{1}{\Gamma\left(m+\dfrac{2-n}{2}\right)}
\hyp32{-n,1,\dfrac{1}{2}}{m+\dfrac{2-n}{2},\dfrac{1-n}{2}}{t^2 }.
\]
Application of \eqref{44.2} gives
\[
k(t)=\dfrac{2^n\Gamma\left(n+\dfrac{3}{2}\right)}{\Gamma\left(\dfrac{1}{2}\right)}
\hyp21{-\dfrac{n}{2},\dfrac{n+1}{2}}{\dfrac{1}{2}}{t^2}.
\]
The hypergeometric function is equal the Legendre polynomial. At last we get
\[
k(t)=(-1)^n\dfrac{\Gamma(2n+2)}{2^{n+1}\Gamma(n+1)}P_n(t).
\]
For $n$ odd we get the same function. This is the expected kernel for the Legendre derivative.

\section{Second method for an alternative for the orthogonal derivative}

Diekema and Koornwinder derive formula for the approximation of higher derivatives in \cite[Section 4]{2}. With a slight adaptation their formula \cite[(4.6)]{2} becomes
\begin{multline*}
h^{-n}\sum_{j=0}^m\left(\int_{-1}^1f(x+ht)P_{n+2j}(t)dt\right)\dfrac{P_{n+2j}^{(n)}(0)}{h_{n+2j}}= \\
=f^{(n)}(x)+(-1)^m\dfrac{\left|P_{n+2m+2}^{(n)}(0)\right|f^{(n+2m+2)}(x)}
{\left|k_{n+2m+2}\right|(n+2m+2)!}h^{2m+2}+o(h^{2m+2}).
\end{multline*}
In this formula $P_j^{(n)}(t)$ denotes the $n$-th derivative of the Legendre polynomial $P_n(t)$. We use the notations $P_n(t)=k_n t^n+\dots$ and $h_n=\int_{-1}^1 P_n(t)^2dt=2/(2n+1)$. We can write the left-hand side as
\[
(-h)^{-n}\int_{-1}^1 f(x+\delta t)k(t)dt
\]
with
\begin{equation}
k(t)=(-1)^n\sum_{j=0}^m \dfrac{P_{n+2j}(t)P^{(n)}_{n+2j}(0)}{h_{n+2j}}.
\label{6.1a}
\end{equation}
With this formula we can derive a simple formula for the function $\omega(t)$. For the Legendre polynomials we have the following formula
\begin{align}
&P_{n+2j}(t)=\dfrac{(-1)^n}{2^n(2n+1)_n}\dfrac{d^n}{dt^n}\left((1-t^2)^nP^{(n,n)}_{2j}(t)\right) \nonumber \\
&\dfrac{d^n}{dt^n}P_{n+2j}(t)=2^{-n}(n+2j+1)_n P^{(n,n)}_{2j}(t).
\label{6.2a}
\end{align}
$P^{(n,n)}_j$ are the Jacobi polynomials with constants $h_j^{(n,n)}$ and $k_j^{(n,n)}$. Using these formula we get
\begin{align*}
k(t)&=2^{-2n}\sum_{j=0}^m \dfrac{(n+2j+1)_n}{(2j+1)_n} \dfrac{1}{h_{n+2j}}
\dfrac{d^n}{dt^n}\left((1-t^2)^nP^{(n,n)}_{2j}(t)P^{(n,n)}_{2j}(0)\right) \\
&=\dfrac{d^n}{dt^n}
\left(2^{-2n}\sum_{j=0}^m \dfrac{(n+2j+1)_n}{(2j+1)_n}
\dfrac{1}{h_{n+2j}}\left((1-t^2)^nP^{(n,n)}_{2j}(t)P^{(n,n)}_{2j}(0)\right)\right).
\end{align*}
Remembering $k(t)=\omega^{(n)}(t)$ gives for $\omega(t)$
\[
\omega(t)=2^{-2n}(1-t^2)^n\sum_{j=0}^m \dfrac{(n+2j+1)_n}{(2j+1)_n}
\dfrac{1}{h_{n+2j}}P^{(n,n)}_{2j}(t)P^{(n,n)}_{2j}(0).
\]
Using
\[
h_{n+2j}\dfrac{2^{2n}(2j+1)_n}{(n+2j+1)_n}=h_{2j}^{(n,n)}
\]
gives
\[
\omega(t)=(1-t^2)^n\sum_{j=0}^m \dfrac{P^{(n,n)}_{2j}(t)P^{(n,n)}_{2j}(0)}{h_{2j}^{(n,n)}}=
(1-t^2)^n\sum_{k=0}^{2m} \dfrac{P^{(n,n)}_k(t)P^{(n,n)}_k(0)}{h_k^{(n,n)}}.
\]
For the summation we can use the formula of Christoffel-Darboux and get
\begin{align*}
\omega(t)&=(1-t^2)^n\dfrac{k_{2m}^{(n,n)}}{h_{2m}^{(n,n)}k_{2m+1}^{(n,n)}}
\dfrac{P_{2m+1}^{(n,n)}(t)P_{2m}^{(n,n)}(0)-P_{2m}^{(n,n)}(t)P_{2m+1}^{(n,n)}(0)}{t-0} \\
&=(1-t^2)^n\dfrac{k_{2m}^{(n,n)}}{h_{2m}^{(n,n)}k_{2m+1}^{(n,n)}}\dfrac{1}{t}
P_{2m+1}^{(n,n)}(t)P_{2m}^{(n,n)}(0).
\end{align*}
The Jacobi polynomials can be written as Gegenbauer polynomials.
\[
\omega(t)=\dfrac{(1-t^2)^n}{t}
\dfrac{k_{2m}^{(n,n)}P_{2m}^{(n,n)}(0)}{h_{2m}^{(n,n)}k_{2m+1}^{(n,n)}}
\dfrac{(n+1)_n}{(2m+n+2)_n}C_{2m+1}^{(n+1/2)}(t).
\]
For the parameters we get
\begin{align*}
&k_{2m}^{(n,n)}=\dfrac{1}{2^{2m}}\binom{4m+2n}{2m} 
\qquad\qquad\qquad\qquad\qquad\qquad
k_{2m+1}^{(n,n)}=\dfrac{1}{2^{2m+1}}\binom{4m+2+2n}{2m+1} \\
&h_{2m}^{(n,n)}=\dfrac{2^{2n+1}}{(4m+2n+1)}\dfrac{\Gamma(2m+n+1)^2}{(2m)!\Gamma(2m+2n+1)} \qquad\quad
P_{2m}^{(n,n)}(0)=\dfrac{(-1)^m}{2^{2m}}\binom{2m+n}{m}.
\end{align*}
Substitution gives after some manipulation with the Pochhammer symbols and the Gamma functions
\begin{align*}
\omega(t)&=\dfrac{(-1)^m}{\pi}\dfrac{\Gamma\left(m+\dfrac{3}{2}\right)\Gamma\left(n+\dfrac{1}{2}\right)}{\Gamma(m+n+1)}\dfrac{(1-t^2)^n}{t}C_{2m+1}^{(n+1/2)}(t) \\
&=\dfrac{(-1)^m}{2^{2m+2n+1}}\dfrac{\Gamma(2m+2)\Gamma(2n+1)}{\Gamma(m+1)\Gamma(n+1)\Gamma(m+n+1)}\dfrac{(1-t^2)^n}{t}C_{2m+1}^{(n+1/2)}(t)
\end{align*}
We found this formula already in section 3. Now we have to compute the basic kernel function $k(t)$. The parameters are
\[
h_{n+2j}=\dfrac{2}{2n+4j+1}
\]
and using \eqref{6.2a}
\begin{equation}
P_{n+2j}^{(n)}(0)=2^{-n}(n+2j+1)_nP_{2j}^{(n,n)}(0)=
(n+2j+1)_n\dfrac{(-1)^j}{2^{2j+n}}\binom{2j+n}{j}.
\label{6.3b}
\end{equation}
Substitution in \eqref{6.1a} gives after some manipulation with the Pochhammer symbols and the Gamma functions
\begin{equation}
k(t)=(-1)^n\dfrac{2^{n-1}}{\sqrt{\pi}}\sum_{j=0}^m (2n+4j+1)\Gamma\left(n+j+\dfrac{1}{2}\right)\dfrac{1}{j!}(-1)^jP_{n+2j}(t).
\label{6.3a}
\end{equation}

\section{The first derivative} 
Looking at formula \eqref{2.10} we see that the hypergeometric function is a polynomial of degree $n$. $n$ is also the order of the derivative. In this section we will discuss the case where $n=1$ so we look for the kernel $k(t)$ for the first derivative.

\

We get from \eqref{2.10} with $n=1$
\begin{equation}
k_m(t)=\dfrac{4}{\pi}\dfrac{\left(m+\dfrac{3}{2}\right)\Gamma\left(m+\dfrac{5}{2}\right)}{\Gamma(m+1)\Gamma(m+2)}
\sum_{j=0}^m \dfrac{(-m)_j\left(m+\dfrac{5}{2}\right)_j}{\left(\dfrac{3}{2}\right)_j }
\dfrac{1}{j!}\big[j-(j+1)t^2\big]t^{2j-1}.
\label{3.1}
\end{equation}
The first few functions (see appendix B), which we call $k_m(t)$ are
\begin{align*}
&k_0(t)=-\dfrac{3}{2}t \\
&k_1(t)=-\dfrac{15}{8}t(5-7t^2) \\
&k_2(t)=-\dfrac{105}{128}t(35-126t^2+99t^4) \\
&k_3(t)=-\dfrac{315}{512}t(105-693t^2+1287t^4-715t^6). \qquad\qquad\qquad\qquad\qquad\qquad\qquad\qquad
\end{align*}
It is easy to check that the function $k_m(t)$ satisfies the following difference equation 
\[
k_{m+1}(t)=\left(1+\dfrac{C_m}{A_m}-\dfrac{1}{A_m}t^2\right)k_m(t)-
\dfrac{C_m}{A_m}k_{m-1}(t)
\]
with
\[
\dfrac{1}{A_m}=\dfrac{(4m+5)(4m+7)}{(2m+2)^2} \qquad \text{   and }\qquad
\dfrac{C_m}{A_m}=\dfrac{(2m+3)^2(4m+7)}{(2m+2)^2(4m+3)}.
\]
Then we may conclude that the function $k_m(t)$ is orthogonal. Looking at equation \eqref{3.1} it turns out that the summation can be performed. It is possible to write the function $k_m(t)$ in a highly simplified form. 

\

We start with \eqref{3.1} and split the right hand side.

\begin{multline*}
k_m(t)=-\dfrac{4}{\pi}\dfrac{\Gamma\left(m+\dfrac{3}{2}\right)\Gamma\left(m+\dfrac{5}{2}\right)}{\Gamma(m+1)\Gamma(m+2)}
\sum_{j=0}^m \dfrac{(-m)_j\left(m+\dfrac{5}{2}\right)_j}{\left(\dfrac{3}{2}\right)_j }
\dfrac{1}{j!}t^{2j+1}+ \\
+(1-t^2)\dfrac{4}{\pi}\dfrac{\left(m+\dfrac{3}{2}\right)\Gamma\left(m+\dfrac{5}{2}\right)}{\Gamma(m+1)\Gamma(m+2)}
\sum_{j=1}^m \dfrac{(-m)_j\left(m+\dfrac{5}{2}\right)_j}{\left(\dfrac{3}{2}\right)_j }
\dfrac{1}{(j-1)!}t^{2j-1}.
\end{multline*}
Changing the summation in the second term we get

\begin{multline}
k_m(t)=-\dfrac{4}{\pi}\dfrac{\Gamma\left(m+\dfrac{3}{2}\right)\Gamma\left(m+\dfrac{5}{2}\right)}{\Gamma(m+1)\Gamma(m+2)}
\sum_{j=0}^m \dfrac{(-m)_j\left(m+\dfrac{5}{2}\right)_j}{\left(\dfrac{3}{2}\right)_j }
\dfrac{1}{j!}t^{2j+1}+ \\
-(1-t^2)\dfrac{8}{3\pi}\dfrac{\left(m+\dfrac{3}{2}\right)\Gamma\left(m+\dfrac{7}{2}\right)}{\Gamma(m)\Gamma(m+2)}
\sum_{j=0}^{m-1}\dfrac{(1-m)_j\left(m+\dfrac{7}{2}\right)_j}{\left(\dfrac{5}{2}\right)_j}\dfrac{1}{j!}t^{2j+1}.
\label{3.2}
\end{multline}
The summations can be done and there appears Jacobi functions.
\[
\sum_{j=0}^m \dfrac{(-m)_j\left(m+\dfrac{5}{2}\right)_j}{\left(\dfrac{3}{2}\right)_j }
\dfrac{1}{j!}t^{2j+1}=t\dfrac{\sqrt{\pi}\Gamma(m+1)}{2\Gamma\left(m+\dfrac{5}{2}\right)}P_m^{(1/2,1)}(1-2t^2).
\]
\[
\sum_{j=0}^{m-1}\dfrac{(1-m)_j\left(m+\dfrac{7}{2}\right)_j}{\left(\dfrac{5}{2}\right)_j}\dfrac{1}{j!}t^{2j+1}=t\dfrac{3\sqrt{\pi}\Gamma(m)}{4\Gamma\left(m+\dfrac{3}{2}\right)}P_{m-1}^{(3/2,2)}(1-2t^2).
\]
Substitution in \eqref{3.2} gives
\begin{equation}
k_m(t)=-\dfrac{2\Gamma\left(m+\dfrac{5}{2}\right)}{\sqrt{\pi}\Gamma(m+2)}
\left[t\, P_m^{(1/2,1)}(1-2t^2)+\left(m+\dfrac{5}{2}\right)(1-t^2)t\, P_{m-1}^{(3/2,2)}(1-2t^2)\right].
\label{3.3}
\end{equation}

To continue there are two possibilities. The first possibility transforms the Jacobi polynomials to Gegenbauer polynomials. These can be transformed into  Legendre polynomials. We get a sum of two Legendre polynomials.

The second possibility uses transform formulas for the Jacobi polynomials and at last we get a single Jacobi polynomial.

\pagebreak

\underline{First possibility.}

The Jacobi polynomials can be transformed  to Gegenbauer polynomials. Using \cite[22.5.25, 22.5.26]{9} gives
\begin{align*}
&t\, P_m^{(1/2,1)}(1-2t^2)=-\dfrac{1}{(2m+3)}C_{2m+1}^{(3/2)}(t) \\
&t\, P_{m-1}^{(3/2,2)}(1-2t^2)=\dfrac{1}{2(2m+3)\left(m+\dfrac{5}{2}\right)}
\dfrac{d}{dt}\left[C_{2m+1}^{(3/2)}(t)\right].
\end{align*}
For the derivative of the Gegenbauer function we use \cite[22.9.2]{9} and get
\[
\dfrac{d}{dt}\left[C_{2m+1}^{(3/2)}(t)\right]=
\dfrac{3}{t}C_{2m}^{(5/2)}(t)-\dfrac{1}{t^2}C_{2m+1}^{(3/2)}(t).
\]
Substitution in \eqref{3.3} gives
\[
k_m(t)=-\dfrac{\Gamma\left(m+\dfrac{3}{2}\right)}{\sqrt{\pi}\Gamma(m+2)}
\left[\left(-1-\dfrac{1-t^2}{2t^2}\right)C_{2m+1}^{(3/2)}(t)+\dfrac{3}{2}\left(\dfrac{1-t^2}{t}\right)C_{2m}^{(5/2)}(t)\right].
\]
For the second Gegenbauer function we use \cite[22.7.3]{9} and get
\[
3(1-t^2)C_{2m}^{(5/2)}(t)=(4+2m)t\, C_{2m+1}^{(3/2)}-(2m+2)C_{2m+2}^{(3/2)}(t).
\]
Substitution gives
\[
k_m(t)=-\dfrac{\Gamma\left(m+\dfrac{3}{2}\right)}{\sqrt{\pi}\Gamma(m+2)}
\left[\left(\dfrac{3}{2}+m-\dfrac{1}{2t^2}\right)C_{2m+1}^{(3/2)}(t)-\dfrac{1}{t}(m+1)C_{2m+2}^{(3/2)}(t)\right].
\]
With \cite[22.7.3]{9}
\[
C_{n-2}^{(\lambda+1)}(t)=\dfrac{(2\lambda+1)t}{2\lambda(1-t^2)}C_{n-1}^{(\lambda)}(t)-
\dfrac{n}{2\lambda(1-t^2)}C_n^{(\lambda)}(t).
\]
we can decrease the upper indexes of the Gegenbauer functions with one. For $\lambda=1/2$ the Gegenbauer polynomial becomes a Legendre polynomial. Using the standard recurrence equation for the Legendre polynomials \cite[22.7.10]{9}
\[
P_{n+1}(t)=\dfrac{2n+1}{n+1}t\, P_n(t)-\dfrac{n}{n+1}P_{n-1}(t)
\]
a lot of times we get at last
\begin{align}
k_m(t)&=\dfrac{1}{B(m+1,m+2)2^{2m+2}}(-1)^{m+1}\, \dfrac{1}{t^2}
\Big[(2m+2)t\, P_{2m+2}(t)+P_{2m+1}(t)\Big] \nonumber \\
&=\dfrac{(2m+3)}{B(m+1,m+2)(4m+5)2^{2m+2}}(-1)^{m+1}\dfrac{1}{t^2}
\Big[(2m+2) P_{2m+3}(t)+(2m+3)P_{2m+1}(t)\Big]
\label{3.4}
\end{align}
where $B(m+1,m+2)$ is the Beta function. We conclude that the kernel function $k_m(t)$ can be written as a sum of two Legendre functions.

\

\underline{Second possibility.}

We repeat here \eqref{3.3}.
\begin{equation}
k_m(t)=-\dfrac{2\Gamma\left(m+\dfrac{5}{2}\right)}{\sqrt{\pi}\Gamma(m+2)}t
\left[P_m^{(1/2,1)}(1-2t^2)+\left(m+\dfrac{5}{2}\right)(1-t^2)
P_{m-1}^{(3/2,2)}(1-2t^2)\right].
\label{3.5}
\end{equation}
We define the form between the square brackets as
\[
K=P_m^{(1/2,1)}(1-2t^2)+\left(m+\dfrac{5}{2}\right)(1-t^2)\, P_{m-1}^{(3/2,2)}(1-2t^2).
\]
For the Jacobi polynomials there are a lot of transformation formula \cite[10.8]{9}. The goal of these transformations is to bring the first order to $3/2$ and the second order to zero. First of all we use \cite[10.8.(33)]{8} a lot of times
\[
2\left(m+\dfrac{\alpha+\beta}{2}+1\right)P_m^{(\alpha,\beta+1)}(1-2t^2)=
\dfrac{(m+\beta+1)}{(1-t^2)}P_m^{(\alpha,\beta)}(1-2t^2)+
\dfrac{(m+1)}{(1-t^2)}P_{m+1}^{(\alpha,\beta)}(1-2t^2)
\] 
and get at last
\begin{align*}
K&=\dfrac{2(m+1)}{(4m+5)}\dfrac{1}{(1-t^2)}
\left(P_m^{(1/2,0)}(1-2t^2)+P_{m+1}^{(1/2,0)}(1-2t^2)\right)+  \\
&+\dfrac{2m(m+1)(2m+5)}{(4m+3)(4m+5)}\dfrac{1}{(1-t^2)}
\left(P_{m-1}^{(3/2,0)}(1-2t^2)+P_m^{(3/2,0)}(1-2t^2)\right)+ \\
&+\dfrac{2m(m+1)(2m+5)}{(4m+5)(4m+7)}\dfrac{1}{(1-t^2)}
\left(P_m^{(3/2,0)}(1-2t^2)+P_{m+1}^{(3/2,0)}(1-2t^2)\right).
\end{align*}
Only the polynomials with first order $1/2$ should be converted into polynomials with order $3/2$. We use \cite[10.8.(35)]{8}
\[
P_m^{(\alpha-1,0)}(1-2t^2)=\dfrac{(m+\alpha)}{(2m+\alpha)}P_m^{(\alpha,0)}(1-2t^2)-
\dfrac{m}{(2m+\alpha)}P_{m-1}^{(\alpha,0)}(1-2t^2).
\]
After gathering the different terms there remains
\begin{align*}
K&=\dfrac{2m(m+1)(2m+3)}{(4m+3)(4m+5)}\dfrac{1}{(1-t^2)}P_{m-1}^{(3/2,0)}(1-2t^2)+ \\
&+\dfrac{2(m+1)(4m^2+10m+3)}{(4m+3)(4m+7)}\dfrac{1}{(1-t^2)}P_m^{(3/2,0)}(1-2t^2)+ \\
&+\dfrac{2(2m+5)(m+1)^2}{(4m+5)(4m+7)}\dfrac{1}{(1-t^2)}P_{m+ 1}^{(3/2,0)}(1-2t^2).
\end{align*}
and there rests all Jacobi polynomials with the same orders. Only the index of the polynomials are different. We use the standard recurrence equation for the Jacobi polynomials \cite[10.8.(11)]{8}
\begin{align*}
&2(m+\alpha)(m+\beta)(2m+\alpha+\beta+2)P_{m-1}^{(\alpha,\beta)}(x)= \\
&=(2m+\alpha+\beta+1)
\left[(2m+\alpha+\beta)(2m+\alpha+\beta+2)x+\alpha^2-\beta^2\right]P_m^{(\alpha,\beta)}(x)- \\
&-2(m+1)(m+\alpha+\beta+1)(2m+\alpha+\beta)P_{m+1}^{(\alpha,\beta)}(x).
\end{align*}
Application gives at last
\[
K=(m+1)P_m^{(3/2,0)}(1-2t^2)=(m+1)(-1)^m P_m^{(0,3/2)}(2t^2-1).
\]	
Then we get for $k_m(t)$
\begin{equation}
k_m(t)=(-1)^{m+1}\dfrac{2\Gamma\left(m+\dfrac{5}{2}\right)}{\sqrt{\pi}\Gamma(m+1)}t\, P_m^{(0,3/2)}(2t^2-1).
\label{3.6}
\end{equation}

\

\begin{remark}
Combining \eqref{3.4} and \eqref{3.6} results in a relation between a special Jacobi polynomial and two Legendre polynomials.
\begin{equation}
P_m^{(0,3/2)}(2t^2-1)=\dfrac{(m+2)2^{m+3}}{(4m+5)}\dfrac{2}{t^3}
\Big[(2m+2) P_{2m+3}(t)+(2m+3)P_{2m+1}(t)\Big].
\label{3.6a}
\end{equation}
\end{remark}

\begin{remark}
For the second order derivative ($n=2$) the formula for $k(t)$ becomes
\[
k_m(t)=(-1)^m\dfrac{4\Gamma\left(m+\dfrac{7}{2}\right)}{\sqrt{\pi}\Gamma(m+1)(4m+7)} \Big[(2m+4)P_m^{(0,3/2)}(2t^2-1)+(2m+3)P_{m+1}^{(0,3/2)}(2t^2-1)\Big].
\]
With \eqref{3.6a} we can express this function $k(t)$ as a summation of Legendre functions of odd orders.
\end{remark}

As we noted earlier, the function $k_m(t)$ for $n=1$ is orthogonal. Here follows an overview of the orthogonality properties of this function.

\

\

\underline{Overview}

- The orthogonal polynomial
\[
k_m(t)=(-1)^{m+1}\dfrac{\Gamma(2m+4)}{\Gamma(m+1)\Gamma(m+2)2^{2m+2}}t\, P_m^{(0,3/2)}(2t^2-1)\qquad\qquad -1 \leq t \leq 1.
\]
- The basic orthogonality property
\[
\int_{-1}^1 k_m(t)k_n(t)\, w(t)dt=h_m \delta_{mn}
\]
\quad with
\[
h_m=\dfrac{\Gamma(2m+4)^2}{2^{4m+3}(4m+5)\Gamma(m+1)^2\Gamma(m+2)^2}.
\]
- The weight function
\[
w(t)=\left(1+\sum_{k=1}^\infty c_k\, t^{2k+1}\right)t^2.
\]
\quad $c_k$ are arbitrary parameters (and can also be omitted).

- The recurrence equation.
\[
k_{m+1}(t)=\left(1+\dfrac{C_m}{A_m}-\dfrac{1}{A_m}t^2\right)k_m(t)- 
\dfrac{C_m}{A_m}k_{m-1}(t)
\]
\quad with
\[
\dfrac{C_m}{A_m}=\dfrac{(2m+3)^2(4m+7)}{(2m+2)^2(4m+3)}  \text{\qquad and \qquad}
\dfrac{1}{A_m}=\dfrac{(4m+5)(4m+7)}{(2m+2)^2}
\]
\quad and
\[
k_{-1}(t)=0  \text{\qquad and \qquad}  k_0(t)=-\dfrac{3t}{2}.
\]
- The differential-difference equation
\[
t(1-t^2)\dfrac{d}{dt}k_m(t)=\left(\dfrac{4m^2+4m+3}{4m+3}-(2m+1)t^2\right)k_m(t)-\dfrac{(2m+3)^2}{4m+3}k_{m-1}(t).
\] 
- The differential equation
\[
t^2(1-t^2)^2\dfrac{d^2}{dt^2}k_m(t)+A_mt(1-t^2)\dfrac{d}{dt}k_m(t)-B_mk_m(t)=0
\]
\quad with
\[
A_m=2(1-2t^2)   \text{\qquad and \qquad} B_m=2(m+2)(2m+1)t^4-2(2m+3)(m+1)t^2+2.
\]

\section{Transfer function of the new derivative}

Liptaj in his paper \cite{4} shows that the formula for the derivatives are much more accurate then the standard Legendre derivatives. He uses some test functions to demonstrate this effect. The derivative can be regarded as a filter. This filter is characterized by his transfer function. The author in his thesis \cite{1} argues that it is recommended to use the transfer function because then you can see 

- the effect of the noise disturbed input signal on the output signal

- the effect of the frequency content of the input signal on the output signal.

In that case you are independent of the test functions. Liptaj uses test functions with very low different frequency contents. In this way it is difficult to decide how accurate the filter is.

In this section we investigate the transfer function of the filter equation \eqref{2.1} of the derivative which will be repeat here
\[
\dfrac{d^n}{dx^n} f(x) =\lim_{h \rightarrow 0}\left(-\dfrac{1}{h}\right)^n\int_{-1}^1 k(t)f(x+h\,t)dt.
\]
Let us call the transfer function of the left hand side $H(\omega)$. The transfer function follows from \eqref{6.3a} and is then given by (omitting the limit)
\begin{multline}
H(h,\omega)=\left(\dfrac{1}{h}\right)^n\dfrac{2^n}{\sqrt{\pi}}
\sum_{j=0}^m (2n+4j+1)\Gamma\left(n+j+\dfrac{1}{2}\right)\dfrac{1}{j!}(-1)^j
\int_{-1}^1 P_{n+2j}(t)e^{-i h\omega t}dt.
\label{4.1}
\end{multline}
The integral is known \cite[3.3(1)]{5}
\begin{align*}
\int_{-1}^1 P_n(t)e^{-i \omega t}dt
&=(-1)^ni^n(2\pi)^{1/2}\omega^{-1/2}J_{n+1/2}(\omega) \\
&=(-1)^ni^n j_n(\omega)
\end{align*}
$J_n(\omega)$ is the Bessel function of the first kind. $j_n(\omega)$ is the spherical Bessel function of the first kind. Substitution gives
\[
H(h,\omega)=\left(-\dfrac{1}{h}\right)^n i^n\dfrac{2^n}{\sqrt{\pi}}
\sum_{j=0}^m (2n+4j+1)\Gamma\left(n+j+\dfrac{1}{2}\right)\dfrac{1}{j!}j_{n+2j}(h\omega)
\]

In the next figure we show the modulus of the transfer function for $n=2$ and different values for $m$. For $m=0$ we get the transfer function of the second order Legendre derivative. This is also given in \cite[section 2.6]{1}. When using another order of the derivative, the shape of the figures does not change, only the axes are adjusted. The same for the value of $h$. The product $(h\omega)$ determines the adjustment of the horizontal axe.
\begin{figure}[ht]
\centering
\includegraphics[height=6cm]{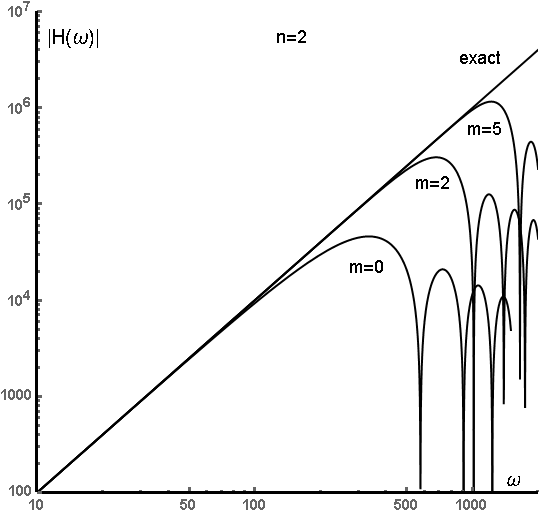}
\caption{Moduli of the absolute value of the transfer function of a second order differentiator for \ $m=0$, \ $m=2$ \ and \ $m=5$. $h=0.01$. For $m=0$ this is the transfer function of the second order Legendre differentiator.}
\end{figure}

We see in the figure that with higher values of $m$ the bandwidth of the derivative will increase. This is in line with what we expect. But there is a also a limit for the bandwidth of the transfer function.

As an approximation of the maximum of the transfer function we get
\[
\omega_{nax} \approx \dfrac{2}{h}\sqrt{2m+n+\dfrac{5}{2}}.
\]
For practical reasons it is recommended using a signal with its highest frquency component is lower then 10 times this maximum. This means that if we want to use test functions, carefully attention must be paid to the frequency content of those test functions. After all the accuracy of the differentiator depends on the frequency content of the test signal. With a test signal with low frequencies the approximation will be better than with a test signal with high frequencies. 

\section{Summary and conclusion} 
This paper describes an extension of the orthogonal derivative which first appears under this name in the thesis of the author \cite{1}. In \cite{2} there is a survey of this derivative which is already known for long time under different names. The standard formula is 
\[
\dfrac{d^n}{dx^n} f(x) =\lim\limits_{h\downarrow 0}\dfrac{k_{n}n!}{h_{n}}
\dfrac{1}{h^{n}}\int_{\RR}f(x+hu) p_{n}(u) d\mu (u)
\]
with the necessary conditions \eqref{I.3}. Liptay in his papers \cite{3},\cite{4} developed a new formula for the derivative by integration of a function on an bounded interval with greater accuracy then the approximate orthogonal derivative. His formula is 
\[
\dfrac{d^n}{dx^n} f(x) =\lim_{h \rightarrow 0}\left(-\dfrac{1}{h}\right)^n\int_{-1}^1 k(t)f(x+h\,t)dt
\]
with $k(t)$ is a kernel which can be computed by solution of a number of equations.
He claims that the accuracy of his formula is much better then the standard Lanczos' derivative which is the first order Legendre derivative. He did not derive a general formula for $k(t)$.
In this paper a general solution for $k(t)$ is derived. We found
\begin{multline*}
k(t)=\dfrac{2}{\pi}\dfrac{\Gamma\left(m+n+\dfrac{3}{2}\right)\Gamma\left(m+\dfrac{3}{2}\right)}{\Gamma(m+1)\Gamma(m+n+1)}
\sum_{j=0}^m \dfrac{(-m)_j\left(m+n+\dfrac{3}{2}\right)_j}{\left(\dfrac{3}{2}\right)_j }
\dfrac{1}{j!}t^{2j-n} \\
\dfrac{\Gamma(2j+1)}{\Gamma(2j+1-n)}
\hyp32{-n,j+1,j+\dfrac{1}{2}}{j+1-\dfrac{n}{2},j+\dfrac{1}{2}-\dfrac{n}{2}}{t^2}.
\end{multline*}
Diekema and Koornwinder  \cite{2} found
\[
k(t)=(-1)^n\dfrac{2^{n-1}}{\sqrt{\pi}}\sum_{j=0}^m (2n+4j+1)\Gamma\left(n+j+\dfrac{1}{2}\right)\dfrac{1}{j!}(-1)^jP_{n+2j}(t).
\]

As a special case we developed a formula for the kernel for the first derivative. We get
\begin{align*}
k_m(t)&=\dfrac{1}{B(m+1,m+2)2^{2m+2}}(-1)^{m+1}\, \dfrac{1}{t^2}\Big[(2m+2)t\, P_{2m+2}(t)+P_{2m+1}(t)\Big] \\
&=(-1)^{m+1}\dfrac{2\Gamma\left(m+\dfrac{5}{2}\right)}{\sqrt{\pi}\Gamma(m+1)}t\, P_m^{(0,3/2)}(2t^2-1).
\end{align*}
In the first formula there are Legendre polynomials and in the second formula there is a Jacobi polynomial. This kernel function is orthogonal. For higher order the kernel function is not orthogonal. For all orders of the derivative it turns out that k$(t)$ can be written as a sum of two Legendre polynomials where each polynomial is multiplied by a polynomial of $t$. Liptaj comes to the same conclusion \cite[Appendix A]{4}. 

\ 

In \cite{1} formula are  developed for the transfer functions of the approximate orthogonal derivative. These are used to see how the accuracy of the function behaves in the frequency domain. In this paper we found
\[
H(h,\omega)=\left(-\dfrac{1}{h}\right)^n i^n\dfrac{2^n}{\sqrt{\pi}}
\sum_{j=0}^m (2n+4j+1)\Gamma\left(n+j+\dfrac{1}{2}\right)\dfrac{1}{j!}j_{n+2j}(h\omega)
\]  

\

Finally, we get the following conclusions:
\begin{enumerate}[\indent {--}]
\setlength\itemsep{0.2em}
\item The formula for the approximate derivatives we use in this paper depends on a parameter $h$. When $h$ goes to zero the exact derivative appears. 
\item The idea of Liptaj for his system of deriving formulas for the derivative by integration is a good idea. Only his method to compute the kernel is not very useful in practice.
\item The formula for the basic kernel in this paper is much more useful to compute this kernel.
\item The formula for the kernel for the first order is orthogonal. The formula for the kernel for higher order derivatives are not.
\item The accuracy of the derivative depends on the first term with the parameter $h$ in the Taylor series of the given function. This term can be chosen freely.
\item To test the accuracy of the formula for the derivative, it is recommended to look at the behaviour of the given function in the frequency domain.
\item As expected the bandwidth of the transfer function increases when the number of terms without the parameter $h$ in the Taylor series of the input function increases.
\end{enumerate}

\

\textbf{Acknowledgements}

The author is very grateful to T.H. Koornwinder and W. Van Assche for their comments on section 6.

\

\begin{appendices}
\addtocontents{toc}{\protect\setcounter{tocdepth}{2}}
\makeatletter
\addtocontents{toc}
{\begingroup
	\let\protect\l@chapter\protect\l@section}
\makeatother
$\text{\bf\Large Appendices}$
\section{An integral with a Gegenbauer polynomial}
In this Appendix we compute the integral
\[
I=\int_{-1}^1 t^{i+2k-n+j}(1-t^2)^{n-j}C^{(n-j+1/2)}_j(t)dt.
\]
The integral is not known, but we can convert this integral into an integral that is known. We use
\begin{equation}
C^{(n-j+1/2)}_j(-t)=(-1)^j C^{(n-j+1/2)}_j(-t).
\label{2.5a}
\end{equation}
Define:
\[
p(t)=t^{i+2k-n+j}(1-t^2)^{n-j}.
\]
Then for the integral we get
\begin{align*}
\int_{-1}^1 p(t)C^{(n-j+1/2)}_j(t)dt
&=\int_{-1}^0 p(t)C^{(n-j+1/2)}_j(t)dt+\int_0^1p(t)C^{(n-j+1/2)}_j(t)dt= \\
&=(-1)^{i+2k-n+j}\int_0^1p(t)C^{(n-j+1/2)}_j(-t)dt+\int_0^1p(t)C^{(n-j+1/2)}_j(t)dt= \\
&=\big[1+(-1)^{i-n}\big]\int_0^1 t^{i+2k-n+j}(1-t^2)^{n-j}C^{(n-j+1/2)}_j(t)dt.
\end{align*}
This integral is known \cite[2.21.1(1)]{6}
\[
\int_0^1x^{\alpha-1}(1-x^2)^{\lambda-1/2}C^{(\lambda)}_{2r}(x)dx=
\dfrac{(-1)^r}{2(2r)!}(2\lambda)_{2r}\Gamma\left(\lambda+\dfrac{1}{2}\right)
\dfrac{\Gamma\left(\dfrac{\alpha}{2}\right)}{\Gamma\left(\dfrac{\alpha+1}{2}+\lambda+r\right)}\left(\dfrac{1-\alpha}{2}\right)_r
\]
with \ \ $\lambda>-\dfrac{1}{2}$  \ \ and \ \ $\alpha>0$. Application of this integral with \ \ $\alpha=i+2k-n+j+1$, \ \ $\lambda=n-j+\dfrac{1}{2}
$ \ \ and $r=\dfrac{j}{2}$ gives

\begin{multline*}
\int_0^1 t^{i+2k-n+j}(1-t^2)^{n-j}C^{(n-j+1/2)}_j(t)dt= \\
=\dfrac{(-1)^{j/2}}{2j!}(2n+1-2j)_j\Gamma(n+1-j)
\dfrac{\Gamma\left(\dfrac{i+2k-n+1+j}{2}\right)}{\Gamma\left(\dfrac{n+i+2k+3}{2}\right)}
\left(\dfrac{n-i-2k-j}{2}\right)_{j/2}.
\end{multline*}
After some manipulations with the Pochhammer symbols and the Gamma functions we obtain
\begin{multline*}
\int_0^1 t^{i+2k-n+j}(1-t^2)^{n-j}C^{(n-j+1/2)}_j(t)dt= \\
=(-1)^j\dfrac{\Gamma(n+1)2^j\left(\dfrac{1}{2}-n\right)_j(i+2m-n+1)_j}{(-2n)_j(i+2k+n+1)j!}
\dfrac{\Gamma\left(\dfrac{i+2k-n+1}{2}\right)}{\Gamma\left(\dfrac{i+2k+n+1}{2}\right)}.
\end{multline*}
Then we get
\begin{multline*}
\int_{-1}^1 t^{i+2k-n+j}(1-t^2)^{n-j}C^{(n-j+1/2)}_j(t)dt= \\
=\left[1+(-1)^{i-n}\right](-1)^j\dfrac{\Gamma(n+1)2^j\left(\dfrac{1}{2}-n\right)_j(i+2m-n+1)_j}{(-2n)_j(i+2k+n+1)j!}
\dfrac{\Gamma\left(\dfrac{i+2k-n+1}{2}\right)}{\Gamma\left(\dfrac{i+2k+n+1}{2}\right)}.
\end{multline*}

\section{Overview of the kernel functions}
In this appendix an overview is given of the kernel functions $k_m(t)$ for different values of the order $n$ and the parameter $m$ using \eqref{2.10}.

\vspace{0.3 cm}
- First order $n=1$
\begin{align*}
&k_0(t)=-\dfrac{3}{2}t \\
&k_1(t)=-\dfrac{15}{8}t(5-7t^2) \\
&k_2(t)=-\dfrac{105}{128}t(35-126t^2+99t^4) \\
&k_3(t)=-\dfrac{315}{512}t(105-693t^2+1287t^4-715t^6) \\
&k_4(t)=-\dfrac{3465}{32768}t(1155-12012t^2+38610t^4-48620t^6+20995t^8) \\
&k_5(t)=-\dfrac{9009}{131072}t(3003-45045t^2+218790t^4-461890t^6+440895t^8-156009t^{10}.\qquad\qquad\qquad
\end{align*}

\vspace{0.3 cm}
- Second order $n=2$
\begin{align*}
&k_0(t)=\dfrac{15}{4}(-1+3t^2) \\
&k_1(t)=-\dfrac{105}{32}(5-42 t^2+45 t^4) \\
&k_2(t)=\dfrac{315}{256}(-35+567 t^2-1485 t^4+1001 t^6) \\
&k_3(t)=-\dfrac{3465}{4096}(105-2772 t^2+12870 t^4-20020 t^6+9945 t^8) \\
&k_4(t)=\dfrac{45045}{65536}(-231+9009 t^2-64350 t^4+170170 t^6-188955 t^8+74613 t^{10})\\
&k_5(t)=-\dfrac{45045}{524288}(3003-162162 t^2+1640925 t^4-6466460 t^6+11904165 t^8-10296594 t^{10}+ \\
&\qquad\qquad\qquad\qquad\qquad\qquad\qquad\qquad\qquad\qquad\qquad\qquad\qquad\qquad\qquad\qquad\ \ \ +3380195 t^{12}).\qquad \ \ \
\end{align*}

\vspace{0.3 cm}
- Third order $n=3$
\begin{align*}
&k_0(t)=-\dfrac{105}{4}t (-3+5 t^2) \\
&k_1(t)=\dfrac{945}{32}t (21-90 t^2+77 t^4) \\
&k_2(t)=-\dfrac{10395}{256} t (-63+495t^2-1001x^4+585x^7) \\
&k_3(t)=\dfrac{45045}{4096} t (693-8580 t^2+30030 t^4-39780 t^6+17765 t^8) \\
&k_4(t)=-\dfrac{135135}{65536} t (-9009+160875 t^2-850850 t^4+1889550 t^6-1865325 t^8+676039 t^{10}) \\
&k_5(t)=\dfrac{2297295}{524288} t (9009-218790 t^2+1616615 t^4-5290740 t^6+8580495 t^8-6760390 t^{10}+ \\
&\qquad\qquad\qquad\qquad\qquad\qquad\qquad\qquad\qquad\qquad\qquad\qquad\qquad\qquad\qquad\qquad +2064825 t^{12}).\qquad \ \ \
\end{align*}

\vspace{0.3 cm}
- Fourth order $n=4$
\begin{align*}
&k_0(t)=\dfrac{945}{16}(3-30 t^2+35 t^4) \\
&k_1(t)=-\dfrac{10395}{64}(-7+135 t^2-385 t^4+273 t^6) \\
&k_2(t)=\dfrac{135135}{2048}(63-1980 t^2+10010 t^4-16380 t^6+8415 t^8) \\
&k_3(t)=-\dfrac{135135}{8192}(-693+32175 t^2-250250 t^4+696150 t^6-799425 t^8+323323 t^{10}) \\
&k_4(t)=\dfrac{2297295}{262144}(3003-193050 t^2+2127125 t^4-8817900 t^6+16787925 t^8-14872858 t^{10}+ \\
&\qquad\qquad\qquad\qquad\qquad\qquad\qquad\qquad\qquad\qquad\qquad\qquad\qquad\qquad\qquad\qquad\ \  +4970875 t^{12}) \ \ \ \\
&k_5(t)=-\dfrac{43648605}{1048576} (-1287+109395 t^2-1616615 t^4+9258795 t^6-25741485 t^8+37182145 t^{10}- \\
&\qquad\qquad\qquad\qquad\qquad\qquad\qquad\qquad\qquad\qquad\qquad\qquad\qquad\qquad\ -26842725 t^{12} +7653825 t^{14}).
\end{align*}

\vspace{0.3 cm}
- Fifth order $n=5$
\begin{align*}
&k_0(t)=-\dfrac{10395}{16}t (15-70 t^2+63 t^4) \\
&k_1(t)=\dfrac{135135}{64} t (-45+385 t^2-819 t^4+495 t^6) \\
&k_2(t)=-\dfrac{675675}{2048}t(1485-20020 t^2+73710 t^4-100980 t^6+46189 t^8) \\
&k_3(t)=\dfrac{11486475}{8192}t (-1287+25025 t^2-139230 t^4+319770 t^6-323323 t^8+119301 t^{10}) \\
&k_4(t)=-\dfrac{218243025}{262144} t (6435-170170 t^2+1322685 t^4-4476780 t^6+7436429 t^8-5965050 t^{10} \\
&\qquad\qquad\qquad\qquad\qquad\qquad\qquad\qquad\qquad\qquad\qquad\qquad\qquad\qquad\qquad\qquad\ \ \  +1847475 t^{12}) \\
&k_5(t)=\dfrac{43648605}{1048576} t (-328185+11316305 t^2-116660817 t^4+540571185 t^6-1301375075 t^8+ \\
&\qquad\qquad\qquad\qquad\qquad\qquad\qquad\qquad\ \ \ +1691091675 t^{10}-1125112275 t^{12}+300540195 t^{14}). \qquad
\end{align*}

\end{appendices}
	
\

\end{document}